\documentclass[smallextended]{svjour3}       
\smartqed  
\usepackage{graphicx}
\usepackage{mathptmx}      
\usepackage[utf8]{inputenc}
\usepackage{amsmath}
\usepackage{amssymb}
\usepackage{graphicx}
\usepackage{hyperref}
\usepackage{xcolor}
\hypersetup{
    colorlinks,
    linkcolor={red!50!black},
    citecolor={blue!50!black},
    urlcolor={blue!80!black}
}
\usepackage{algorithm}
\usepackage{algorithmic}
\usepackage{subcaption}
\usepackage{tikz}
\usetikzlibrary{shapes,decorations,arrows,calc,arrows.meta,fit,positioning}
\tikzset{
 -Latex,auto,node distance =1 cm and 1 cm,semithick,
 state/.style ={ellipse, draw, minimum width = 0.7 cm},
 point/.style = {circle, draw, inner sep=0.04cm,fill,node contents={}},
 bidirected/.style={Latex-Latex,dashed},
 el/.style = {inner sep=2pt, align=left, sloped},
}
\ifdefined\theorem\else \newtheorem{theorem}{Theorem}\fi
\ifdefined\proposition\else \newtheorem{proposition}[theorem]{Proposition}\fi
\ifdefined\definition\else \newtheorem{definition}[theorem]{Definition}\fi
\ifdefined\lemma\else \newtheorem{lemma}[theorem]{Lemma}\fi
\ifdefined\corollary\else \newtheorem{corollary}[theorem]{Corollary}\fi
\ifdefined\remark\else \newtheorem{remark}[theorem]{Remark}\fi
\ifdefined\example\else \newtheorem{example}{Example}\fi

\newcommand{\R}{\mathbb{R}}

\newcommand{\Z}{\mathbb{Z}}
\newcommand{\N}{\mathbb{N}}
\newcommand{\E}{\mathbb{E}}

\newcommand{\X}{\mathbb{X}}

\newcommand{\Pro}{\mathbb{P}}
\newcommand{\ce}[1]{[\![#1]\!]}         
\newcommand{\np}[1]{\left(#1\right)}
\newcommand{\nc}[1]{\left[#1\right]}

\newcommand{\Y}{\mathbb{Y}}

\newcommand{\argmin}{arg\,min}

\begin{document}

\title{Entropic regularization of the Nested Distance
}


\author{Zheng Qu        \and
 Beno\^it Tran 
}


\institute{Zheng Qu \at
 Department of Mathematics,
 The University of Hong Kong
 Room 419, Run Run Shaw Building
 Pokfulam Road, Hong Kong \\
 \email{zhengqu@hku.hk}           
 \and
 Beno\^it Tran \at
 \'Ecole des Ponts ParisTech, CERMICS,
 6 et 8 avenue Blaise Pascal
 Cit\'e Descartes - Champs-sur-Marne
 77455 Marne-la-Vall\'ee C\'edex 2\\
 \email{benoit.tran@enpc.fr}
}

\date{Received: date / Accepted: date}

\maketitle

\begin{abstract}
 In 2012, Pflug and Pichler proved, under regularity assumptions, that the value function in Multistage Stochastic Programming (MSP) is Lipschitz continuous w.r.t. the Nested Distance, which is a distance between scenario trees (or discrete time stochastic processes with finite support). The Nested Distance is a refinement of the Wasserstein distance to account for proximity of the filtrations of discrete time stochastic processes.
 
 The computation of the Nested Distance between two scenario trees amounts to the computation of an exponential (in the horizon $T$) number of optimal transport problems between smaller conditional probabilities of size $n$, where $n$ is less than maximal number of children of each node.
 
 Such optimal transport problems can be solved by the auction algorithm with complexity $O\np{n^3\log(n)}$. In 2013, Cuturi introduced Sinkhorn's algorithm, an alternating projection scheme which solves an entropic regularized optimal transport problem. Sinkhorn's algorithm converges linearly and each iteration has a complexity of $O(n^2)$.
 
 In this article, we present and test numerically an entropic regularization of the Nested Distance.
 \keywords{Nested Distance \and Sinkhorn's Algorithm \and Multistage Stochastic Programming \and Optimal Transport}
\end{abstract}


\section{Introduction: from the Wasserstein distance to the Nested Distance}

In Multistage Stochastic Programming (MSP), Georg Pflug introduced in 2009 \cite{Pf2009} the \emph{Nested Distance}, which is a refinement of the Wasserstein distance to account proximity in the filtrations between two discrete time stochastic processes. Following usual denomination in the Stochastic Programming community (see \cite{He.Ro2009,Pf.Pi2014,Sh.De.Ru2009}), we also denote by \emph{scenario tree} a discrete time stochastic process which is also discrete and finite in space.

There are many different distances between scenario trees. However, few are suited for MSP purposes: one would like to guarantee continuity of the value function of a MSP with respect to scenario trees, \emph{i.e.} if two scenario trees are arbitrarily close to each other, then the value of the associated MSP (with the same structure except for the scenario trees) can be made arbitrarily close as well.

One possible distance between scenario tree is the Wasserstein distance. Intuitively, the Wasserstein distance between two probabilities $p$ and $q$ (for scenario tree $\np{X_t}_{t\in \ce{1,T}}$, consider the probability law of the tuple $\np{X_1, \ldots, X_T}$) corresponds to the optimal cost of splitting and transporting the mass from one to the other. We write $\mathbf{1}_k$, $k\in \N$, for the vector $\np{1,\ldots, 1}^T$ of $\R^k$.

\begin{definition}[Discrete optimal transport and Wasserstein distances]
 \label{def:OT}
 Let $n,m$ be two integers and $\X = \left\{ x_1, x_2, \ldots, x_n\right\}$ and $\Y = \left\{ y_1, \ldots, y_m\right\}$ be two finite sets included in $\R^N$, $N\geq 1$. Denote by $c=\np{c_{ij}}_{i,j}$ a $n\times m$ positive matrix called \emph{cost matrix}. The \emph{optimal transport cost} between two probability measures $p$ and $q$ on respectively $\X$ and $\Y$, is the value of the following  optimization problem
 \begin{equation}
  \label{eq:OT}
  \mathrm{OT}\np{p, q; c} = \min_{\pi \in \R_+^{n\times m}}    \sum_{\substack{1\leq i \leq n \\ 1\leq j \leq m}} c_{ij}\pi_{ij}  \text{ s.t. }  \pi\mathbf{1}_m = p \text{ and } \pi^T\mathbf{1}_n = q.
 \end{equation}
 Moreover, defining the \emph{cost function} by $c\np{x_i, y_j} = c_{ij}$ for every indexes $i,j$, Problem~\eqref{eq:OT} can be written using probabilistic vocabulary as
 \begin{equation}
  \label{OT_Pro}
  \mathrm{OT}\np{p, q; c} = \min_{\substack{\np{X,Y} \ \text{s.t.}\\X \sim p \\ Y \sim q}} \E_{\np{X,Y}} \nc{c\np{X,Y}},
 \end{equation}
 where the notation $X \sim p$ (resp. $Y \sim q$) means that the probability law of the random variable $X\in \X$ (resp. $Y \in \Y$) is equal to $p$ (resp. $q$) and the notation $\E_{\np{X,Y}}$ is the expectation under the probability law of the couple of random variables $\np{X,Y}$.
 
 Lastly, when for some real $r\geq 1$, the cost function $c$ is equal to $d^r$ with $d$ a metric on $\R^N$, then $\mathrm{OT}\np{p,q; d^r}^{1/r}$ is the \emph{$r$-th Wasserstein distance} between $p$ and $q$, denoted $\mathrm{W_r}\np{p, q}$.
\end{definition}
We refer to the textbooks \cite{Pe.Cu2019,Vi2009} for a presentation and references on optimal transport. In two stage multistage optimization problems, under some regularity assumptions, the value function is Lipschitz continuous with respect to the Wasserstein distances, see \cite[Chapter 6]{Pf.Pi2014}. However the value function of MSP with more than $2$ stages is not continuous with respect to the Wasserstein distances, as seen in Example~\ref{Wasserstein_and_MSP}, where we show that for a $3$ stage MSP, two scenario trees can be arbitrarily close to each other in the $1$-Wasserstein metric but the gap in value of the associated MSP is arbitrarily large.

\begin{example}[The Wasserstein distance is not suited for MSP]
 \label{Wasserstein_and_MSP}
 In this example we illustrate that the $1$-Wasserstein is not an interesting metric to evaluate distance between scenario trees involved in a MSP: an arbitrary small Wasserstein distance between two scenario trees may yield an arbitrary large gap in values of the same MSP.
 
 Given a scenario tree $Z$ (see Definition~\ref{def:scenariotree} for a formal definition) with natural filtration $\np{\mathcal{F}_t}_{t\in \ce{0,2}}$\footnote{For every $t\in \ce{0,2}$, $\mathcal{F}_t = \sigma\np{Z_0,\ldots, Z_t}$.}, we want to buy a single object at the minimal average cost
 \[
  v\np{Z} = \min_{\mathbf{u}} \left\{\mathbb{E}\left[\sum_{t=0}^{2} Z_{t} \mathbf{u_t}\right] \mid \begin{array}{c}{ \mathbf{u_t} \in \left\{ 0, 1\right\},} \\ {\mathbf{u_t} \text { is } \mathcal{F}_{t} \text { -measurable, }}  \\ \sum_{t=0}^T \mathbf{u_t} = 1, \end{array}\right\}.
 \]
 
 Fix $A \gg \epsilon > 0$, in Figure~\ref{fig:WassersteinNotSuited} are two scenario tree modeling the price of an object during $3$ time steps. Their natural filtrations are different. Intuitively, on the left scenario tree, the decision maker observes that an $\epsilon$ variation of the price at $t=1$ and knows that it will yield an explosion (upward or downward) of the price at $t=2$. Whereas on the right scenario tree, the decision maker does not recognize such information at time $t=1$.  Example inspired from \cite{He.Ro.St2006}.
 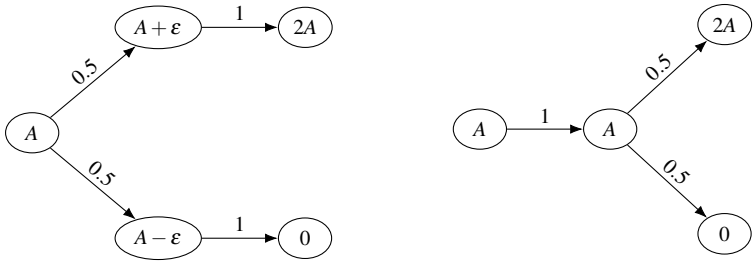
\begin{figure}[ht]
  \begin{minipage}[c]{.49\linewidth}
   \begin{tikzpicture}
    \node[state] (1) {$A$};
    \node[state] (2) [above right =of 1] {$A + \epsilon$};
    \node[state] (3) [below right=of 1] {$A - \epsilon$};
    \node[state] (4) [right =of 2] {$2A$};
    \node[state] (5) [right =of 3] {$0$};
    \path (1) edge node[el,above] {$0.5$} (2);
    \path (1) edge node[el,above] {$0.5$} (3);
    \path (2) edge node[above] {$1$} (4);
    \path (3) edge node[above] {$1$} (5);
   \end{tikzpicture}
  \end{minipage}
  \begin{minipage}[c]{.5\linewidth}
   \begin{tikzpicture}
    \node[state] (1) {$A$};
    \node[state] (2) [right =of 1] {$A$};
    \node[state] (3) [above right =of 2] {$2A$};
    \node[state] (4) [below right =of 2] {$0$};
    \path (1) edge node[el, above] {$1$} (2);
    \path (2) edge node[el, above] {$0.5$} (3);
    \path (2) edge node[el, above] {$0.5$} (4);
   \end{tikzpicture}
  \end{minipage}
  \caption{\label{fig:WassersteinNotSuited}Left: scenario tree $X := \np{X_0, X_1, X_2}$. Right: scenario tree $Y = \np{Y_0, Y_1, Y_2}$.}
 \end{figure}
 
 On the one hand we have proximity in the $1$-Wasserstein metric $\mathrm{W}_1$ as
 \(
 \mathrm{W}_1\np{X, Y} = 2\epsilon.
 \)
 On the other hand, the optimal values are $v\np{X} = \frac{A+\epsilon}{2}$ and $v\np{Y} = A$. Thus, we have an arbitrarily large gap in values
 \[
  \lvert v\np{X} - v\np{Y} \rvert = \frac{A-\epsilon}{2} \underset{A\to +\infty}{\rightarrow} +\infty.
 \]
\end{example}

In 2012, Pflug and Pichler proved in \cite{Pf.Pi2012} that the Nested Distance previously introduced by Pflug, is the correct adaptation of the Wasserstein distance for multistage stochastic programming: under regularity assumptions, the value function of MSPs is Lipschitz continuous with respect to the Nested Distance between scenario trees. Since then, it has been used as a tool to quantify the quality of approximating trees: given an initial scenario tree, one would like to have a good approximating tree with fewer nodes. The Nested Distance both quantifies the quality of an approximating tree and the associated optimal transport plan also allows for reduction of scenario trees,  see for example \cite{Ho.Vi.Ko.Mo2020,Ko.Pi2015}.

Without additional structure (like independence) of the scenario tree, the Nested Distance is usually computed via a backward recursive algorithm (introduced in~\cite{Pf.Pi2012}, see also \cite[Definition 15]{Pi.Sc2019}) which amounts to solve an exponential number (in $T$) number of optimal transportation problems. It decomposes over the time the computation of the Nested Distance as the dynamic computation of a finite number of optimal transport problems between conditional probabilities with costs updated backward.

Optimal transport between discrete probabilities of size $n$ can be solved by Linear Programming algorithms like the auction algorithm with complexity roughly $O\np{n^3\log{n}}$, see \cite{Be.Ca1989}.

\textbf{Main idea.} By adding an entropic term to the primal of the optimal transport problem associated with the computation of a Wasserstein cost, an alternating projection scheme yield Sinkhorn's algorithm, introduced in Optimal Transport in \cite{Cu2013} to compute Wasserstein distances. The complexity of each iteration of Sinkhorn's algorithm is $O\np{n^2}$ and Sinkhorn's algorithm generates a sequence of transport plans which converges linearly to the minimizer of the regularized Optimal Transport problem. Relaxing each optimal transport problem involved in the recursive computation of the Nested Distance, we end up with an entropic regularization of the Nested Distance.

The remainder of the article is organized as follows:
\begin{itemize}
 \item In Section~\ref{sec:END}, we first formally define the Nested Distance as the value of a dynamic system of optimal transport problems between conditional probabilites and varying costs. Then, we present an entropic regularization of the discrete optimal transport Problem~\eqref{eq:OT} and how this relaxed OT problem can be solved efficiently by Sinkhorn's algorithm. Lastly, we define a natural entropic regularization of the Nested Distance by relaxing each OT problem involved in its dynamic formulation.
 \item In Section~\ref{sec:num}, we end this article with a numerical experiment showing both the speedup of our approach to compute Nested Distances and also its relative preciseness.
\end{itemize}

\textbf{Related work.} The present work has first been presented in one chapter of the second author's Ph.D. thesis \cite{tran:tel-03129146}. Since then, Pichler and Weinhardt in \cite{Pi.We2021} extended this work by deriving a dual characterization of the regularized Nested Distance and also gave an upper bound on the approximation error when approximating the Nested Distance by its regularized counterpart for a given regularization parameter.

\begin{figure}
 \begin{subfigure}[b]{0.49\textwidth}
  \includegraphics[width=\linewidth]{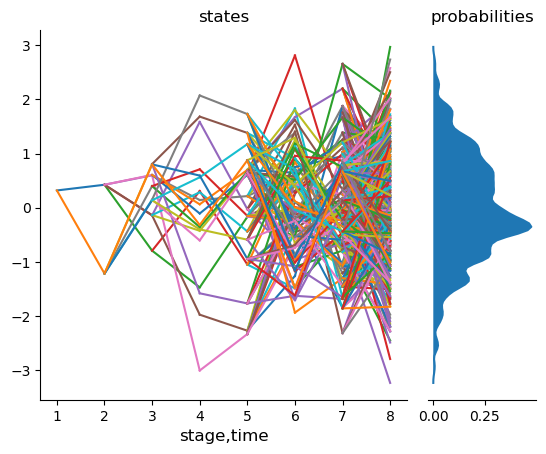}
 \end{subfigure}
 \begin{subfigure}[b]{0.49\textwidth}
  \includegraphics[width=\linewidth]{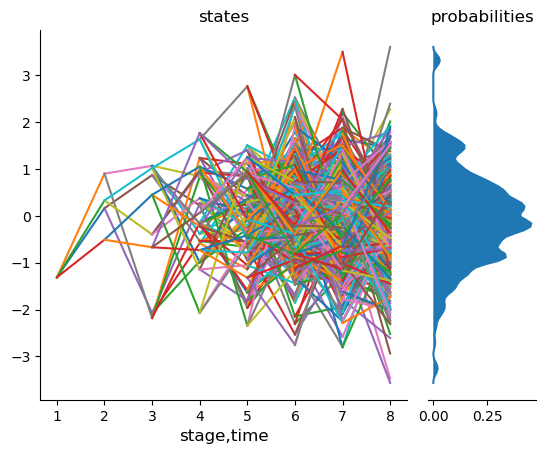}
 \end{subfigure}
 \caption{Two scenario trees $X$ and $Y$ with a continuous probability approximation of the histogram the leaves. Their Nested Distance is $\mathrm{ND}_2\np{X, Y} = 1.009$ and its entropic regularization is $\mathrm{END}_2\np{X,Y} = 1.011$, see Section~\ref{sec:num}. The trees were generated using the ScenTrees.jl package \cite{Ki.Pi.Pf2020}.}
\end{figure}

\section{The Nested Distance and its entropic regularization}
\label{sec:END}

\subsection{Dynamic computation of the Nested Distance}
Throughout the remainder of the article, we fix an integer $T > 1$ and we consider that the set of reals $\R$ is endowed with its usual distance and borelian structure. Moreover, for every $t\in \ce{1,T}$, $\R^t = \R \times \ldots \times \R$ is seen as a filtered space endowed with its cylinder $\sigma$-algebra.

\begin{definition}[Scenario tree]
 \label{def:scenariotree}
 Let $\np{X_t}_{t\in \ce{1,T}}$ be a discrete time stochastic process defined on some probability space. The stochastic process $\np{X_t}_{t\in \ce{1,T}}$ is a \emph{scenario tree} if it is also finite and discrete in space, \emph{i.e.} for every time indexes $1\leq s\leq t\leq T$, the support $\X_{s:t}$ of $X_{s:t} = \np{X_s,\ldots, X_t}$ defined by
 \[
  \X_{s:t} := \left\{ x_{s:t} = \np{x_s,\ldots, x_t} \in \R^{t-s} \mid  \Pro\np{X_s = x_s, \ldots, X_t = x_t} > 0\right\}
 \]
 is non-empty, finite and
 \(
 \sum_{x_{s:t} \in \X_{s:t}} \Pro\np{X_s = x_s, \ldots, X_t = x_t} = 1.
 \)
\end{definition}

Following \cite{Pi.Sc2019}, we define the Nested Distance between scenario trees as the value of a recursive computation of optimal transportation between conditional probabilities with updated costs. Given two scenario trees $X = \np{X_t}_{t\in \ce{1,T}}$ and $Y =\np{Y_t}_{t\in \ce{1,T}}$, for every $s,t\in \ce{1,T}$, we define the tuple of random variable variables $X_{s:t} = \np{X_s, \ldots X_t}$ and $Y_{s:t} = \np{Y_s, \ldots, Y_t}$. Denote by $x_{s:t}$ and $y_{s:t}$ any element of their support $\X_{s:t}$ and $\Y_{s:t}$, defined in Definition~\ref{def:scenariotree}. Lastly, for every $t\in\ce{1,T}$, denote by $P_t$ and $\tilde{P}_t$ the probability law of $X_{1:t} = \np{X_1, \ldots, X_t}$ and $Y_{1:t} = \np{Y_1, \ldots, Y_t}$, respectively.
\begin{definition}[Nested Distance between scenario trees]
 \label{def:ND}
 Let $X$ and $Y$ be two scenario trees. Given $r\geq 1$, and the metric $d\np{x,y} = \Vert x-y \Vert_r$ over $\R^T$, for every $t\in \ce{1,T}$, compute recursively backward in time functions $c_t : \X_{1:T} \times \Y_{1:T} \to \overline{\R}$ by
 \begin{equation}
  \label{NestedDistance}
  \left\{
  \begin{aligned}
    & c_T\np{x_{1:T}, y_{1:T}} = d\np{x_{1:T}, y_{1:T}}, \ \forall \np{x_{1:T}, y_{1:T}} \in \X_{1:T} \times \Y_{1:T},                                             \\
    & c_{t}\np{x_{1:T}, y_{1:T}} = \mathrm{OT}\np{P_{t+1}\np{ \cdot \mid X_{1:t} = x_{1:t}}, \tilde{P}_{t+1}\np{\cdot \mid Y_{1:t} = y_{1:t}}; c_{t+1}^r } ^{1/r}, \\
    & \forall t\in \ce{1,T-1}, \ \forall \np{x_{1:T}, y_{1:T}} \in \X_{1:T} \times \Y_{1:T}.
  \end{aligned}
  \right.
 \end{equation}
 Set $\mathrm{ND}_r\np{X, Y} := \mathrm{OT}\np{P_T, \tilde{P}_T, c_1^r}^{1/r}$, it is the \emph{$r$-Nested Distance} between the scenario trees $X$ and $Y$.
\end{definition}

Although for every $t\in \ce{1,T}$ the domain of $c_t$ is $\X_{1:T}\times \Y_{1:T}$, only the process up to $t$ matters \emph{i.e.} for every $x_{1:T},x'_{1:T} \in \X_{1:T}$ and $y_{1:T},y'_{1:T} \in \Y_{1:T}$ such that $x_{1:t} = x'_{1:t}$ and $y_{1:t} = y'_{1:t}$ we have $ c_{t}\np{x_{1:T}, y_{1:T}} = c_{t}\np{x'_{1:T}, y'_{1:T}}$.
It follows from \cite[Proposition 20]{Pi.Sc2019} and \cite[Theorem 19]{Pf.Pi2012} that the Nested Distances $\mathrm{ND}_r$ introduced in Definition~\ref{def:ND}, are distances on the space of scenario trees.

\begin{remark}
 Solving Problem~\eqref{NestedDistance} amounts to solving an exponential (in $T$) number of Linear optimization Problems where the dimension of the variable to optimize is bounded by $\max_{t\in \ce{1:T-1}} \lvert {\X_{t:t+1}} \rvert \cdot \max_{t\in \ce{1:T-1}} \lvert \Y_{t:t+1} \rvert $.
\end{remark}

\subsection{Entropic regularization of optimal transport problems}
\label{entropicOT}
We will regularize the OT Problem~\ref{eq:OT} by adding an entropy term to the objective function. The \emph{Shannon entropy} or simply \emph{entropy} of a random variable $Z$ with values in a finite subset $\Z$ of cardinal $k\in \N$ in $\R^t$, $t\geq 1$ and probability vector $\np{p_1;\ldots; p_k} \in \np{\R_{+}^*}^k$ is defined as
\(
H(Z) = \E\nc{ - \log Z} = - \sum_{i =1}^k p_i \log\np{p_i}.
\)
By adding an entropy regularization term to the objective of an optimal transport Problem~\ref{OT_Pro} (using the probabilistic notations), the linear objective function of a discrete OT problem becomes strongly convex, hence damping the combinatorial aspects of OT.

\begin{definition}[Regularized Optimal Transport]
 With the notations of Definition~\ref{def:OT}, for every real $\gamma > 0$ we define the following \emph{regularized optimal transport plan} between probabilities $p \in \R^n$ and $q \in \R^m$ with cost matrix $c \in \R^{n\times m}$
 \begin{equation}
  \label{eq:reg_OT}
  \pi_{\gamma}\np{p,q;c} = \argmin_{\substack{\np{X,Y} s.t. \\ X \sim p \\ Y \sim q}} \E\nc{c\np{X,Y} - \gamma H(\mathcal{L}\np{X,Y})},
 \end{equation}
 where $\mathcal{L}\np{X,Y}$ is the probability law of the couple $\np{X,Y}$ of random variables. Then, the associated value is the \emph{regularized optimal transport} $\mathrm{OT}_{\gamma}$ between $p$ and $q$
 \begin{equation}
  \mathrm{OT}_{\gamma}\np{p,q;c} = \sum_{\substack{1\leq i \leq n \\ 1\leq j \leq m}} c_{ij}\np{\pi_{\gamma}}_{ij}.
 \end{equation}
\end{definition}
\begin{remark}
 \label{rem:ot_vs_regot}
 Note that as the regularized optimal transport plan $\pi_{\gamma}$ also satisfies the constraints of the (unregularized) optimal transport problem, we have for every $\gamma > 0$ that $\mathrm{OT}\np{p,q;c} \leq \mathrm{OT}_{\gamma}\np{p,q;c}$. Moreover when $\gamma$ tends to $0$, one recover the optimal transport value, \emph{i.e.}
 \(
 \mathrm{OT}_{\gamma}\np{p,q;c} \underset{\gamma \to 0}{\longrightarrow} \mathrm{OT}\np{p,q;c}.
 \)
\end{remark}
Given an integer $t\geq 1$, let $p$ and $q$ be two probabilities on $\R^t$ with respective finite support of size $n\in \N$ and $m \in \N$. We say that a $n\times m$ matrix $\pi$ is a  \emph{transport plan between $p$ and $q$} if it is admissible in Problem~\eqref{eq:OT}, \emph{i.e.} $\pi$ satisfies the \emph{mass conservation constraints}:
\begin{equation}
 \label{eq:mass_conservation}
 \pi\mathbf{1}_m = p \text{  and  }  \pi^T\mathbf{1}_n = q.
\end{equation}
The set of transport plans between $p$ and $q$ is denoted by $\mathcal{P}\np{p,q}$. We now present Sinkhorn's algorithm. This algorithm was (re)discovered by Cuturi in \cite{Cu2013} who used it to solve the regularized optimization Problem~\eqref{eq:reg_OT}. Proofs of the different following statements can be found in \cite{Cu2013}, \cite[Chapter 4]{Pe.Cu2019} and \cite{Al.Ni.Ri2017}.

\begin{theorem}[Sinkhorn's algorithm and its convergence]
 \label{theo:Sinkhorn}
 Fix $\gamma > 0$, an integer $t\geq 1$ and let $p$ and $q$ be two probabilities on $\R^t$ with respective finite support of size $n\in \N$ and $m \in \N$. The following assertions are true:
 \begin{itemize}
  \item  There exists a unique transport plan $\pi^*$ which minimizes the regularized optimal transport Problem~\eqref{eq:reg_OT} with cost matrix $c=\np{c_{ij}}_{i,j} \in \R^{n\times m}$.
  \item  There exists vectors $u^*\in \np{\R_+^*}^n$, $v^* \in \np{\R_+^*}^m$ such that
        \(
        \pi^* = \mathrm{diag}\np{u^*}G \,\mathrm{diag}\np{v^*},
        \)
        where $G$ is the Gibbs kernel defined by $G_{ij} = \exp\np{-\frac{c_{ij}}{\gamma}}$.
  \item Alternatively rescaling the lines and columns of $G$ in order to satisfy the mass conservation constraints of Equation~\eqref{eq:mass_conservation} converges to the optimal transport plan $\pi^*$. More precisely, iterates $\np{u_k, v_k} \in \np{\R_+^*}^n \times \np{\R_+^*}^m$, $k\in \N$, defined by $u_0 = \mathbf{1}_n$, $v_0 = \mathbf{1}_m$, and
        \begin{equation}
         \label{eq:Sinkhorn}
         \left\{
         \begin{aligned}
           & u_{k+1} = \mathbf{1}_{n} \ ./ \np{Gv_k}   \quad \text{(where $./$ is the entrywise division)} \\
           & v_{k+1} = \mathbf{1}_m\ ./ \np{G u_{k+1}},
         \end{aligned}
         \right.
        \end{equation}
        converge  to the optimal scaling vectors $u^*$ and $v^*$.
  \item Sinkhorn's iterates $\np{u_k, v_k}_{k\in \N}$ defined in Equation~\eqref{eq:Sinkhorn} converge linearly to the optimal scaling vectors $\np{u^*, v^*}$.
 \end{itemize}
\end{theorem}

\begin{figure}[ht]
 \leftline{\begin{subfigure}[b]{0.33\textwidth}
   \includegraphics[width=\linewidth]{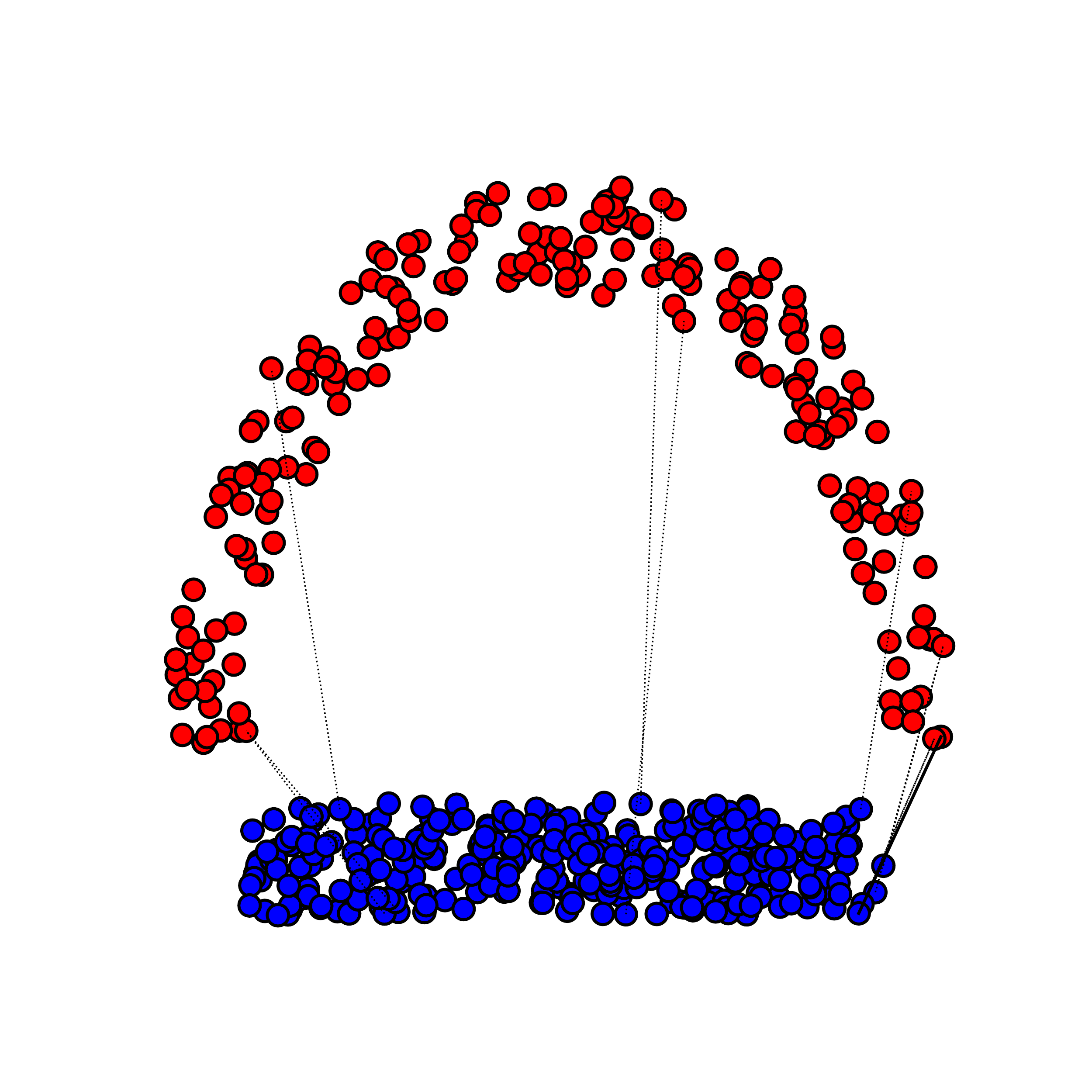}
  \end{subfigure} \begin{subfigure}[b]{0.33\textwidth}
   \includegraphics[width=\linewidth]{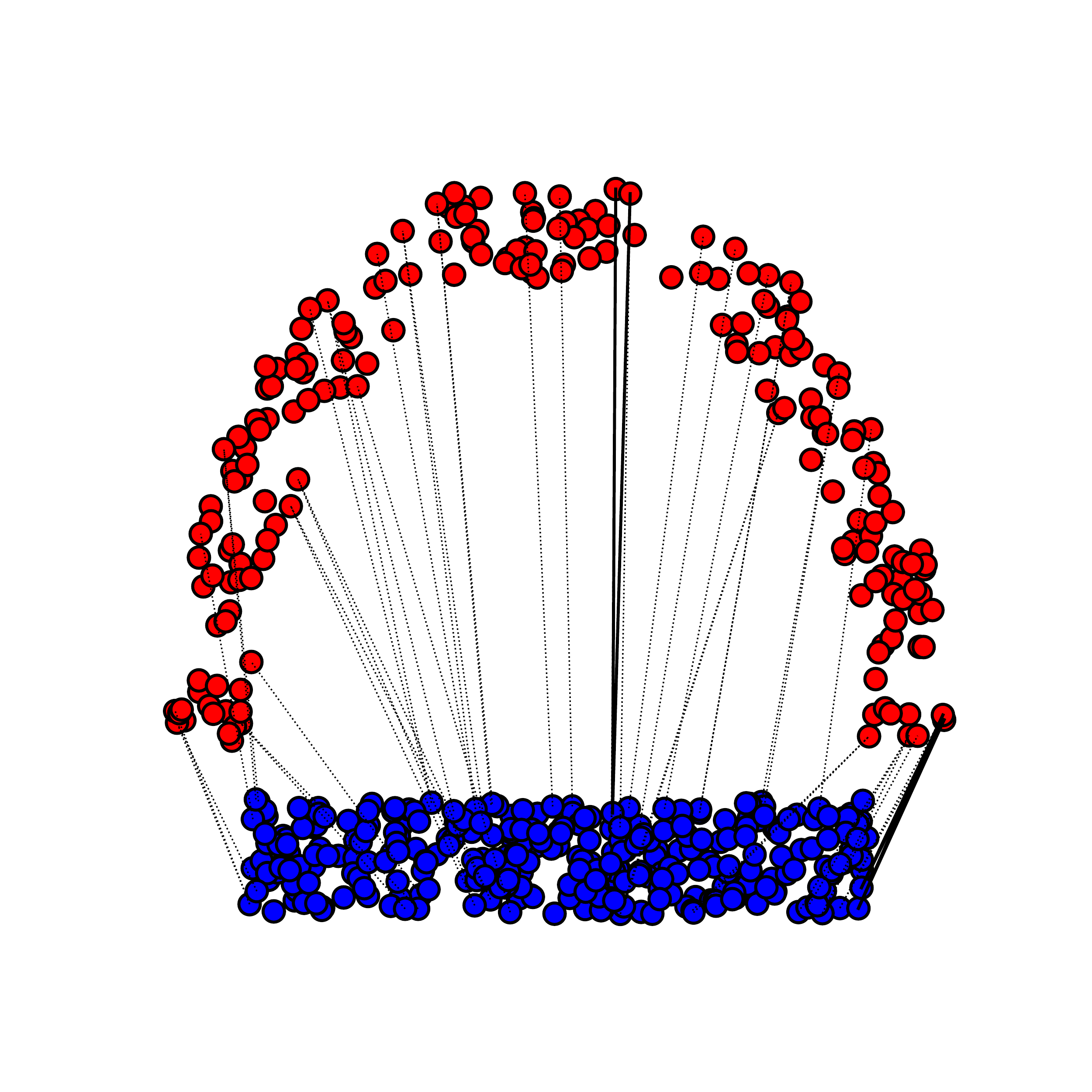}
  \end{subfigure}
  \begin{subfigure}[b]{0.33\textwidth}
   \includegraphics[width=\linewidth]{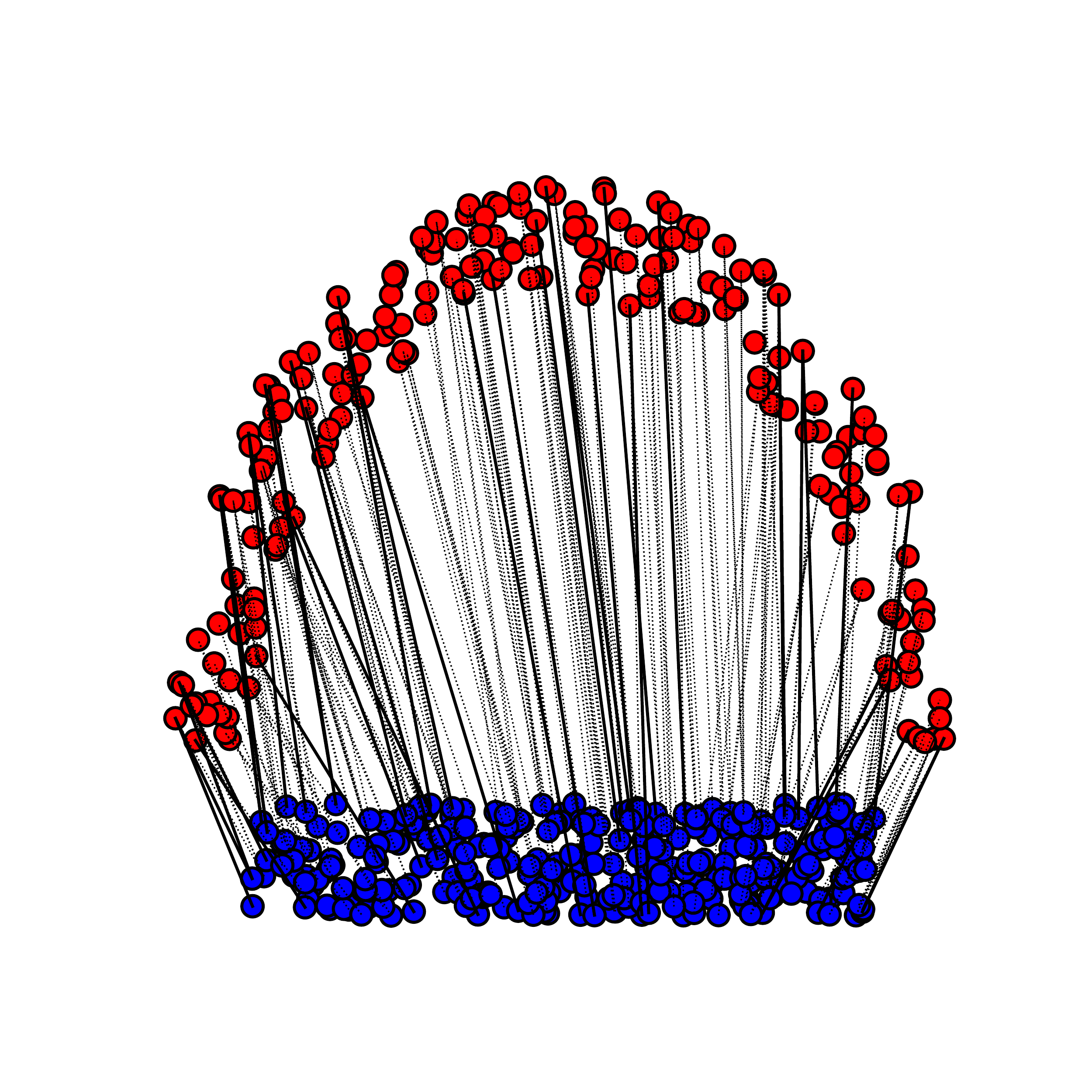}
  \end{subfigure}}%
 \caption{\label{fig:opt_relaxedtransportplan}Left: $\gamma = 0.008$. Middle: $\gamma = 0.005$. Right: $\gamma = 0.003$.  Effect of the regularization parameter on the optimal transport plan of the relaxed OT Problem~\eqref{eq:reg_OT}, between the empirical probability measures of the red cloud and the blue cloud. In all cases, a continuous edge (resp. dashed line) exists if more than $30$ percent (resp. $20$ percent) of the red dot mass is moved to the associated blue one.}
\end{figure}

We comment on the optimal transport plan associated with the relaxed OT Problem~\eqref{eq:reg_OT}. When the regularization parameter $\gamma$ is large, then the optimal transport plan is very diffuse: in the left part of Figure~\ref{fig:opt_relaxedtransportplan}, this means that the mass of each red dot is spread along many different blue dots. The closer the regularization parameter $\gamma$ gets to $0$, then the combinatorial aspect of discrete OT appears gradually: each red dot is spread along few different blue dots. This is expected, as if there were the same number of blue and red dots in Figure~\ref{fig:opt_relaxedtransportplan}, then the OT problem is an assignment problem. Informally, the entropic regularization of Equation~\eqref{eq:reg_OT} dampens the combinatorial aspects of the optimal transport problem of Equation~\eqref{eq:OT}.

\subsection{Entropic regularization of the Nested Distance}

We have seen in \S\ref{entropicOT} how to compute efficient upper bound $\mathrm{OT}_{\gamma}$ of the discrete optimal transport problem $\mathrm{OT}$. Hence, by replacing optimal transport problems by their relaxed counterpart in the dynamic computation of the Nested Distance in Equation~\eqref{NestedDistance}, we have an entropic regularization of the Nested Distance noted $\mathrm{END}$.

Note that a regularizing parameter $\gamma > 0$ must be chosen for each optimal transport problem in Equation~\eqref{NestedDistance}. On the one hand, one would like to put $\gamma$ as small as possible in order to have the best approximation of the unregularized OT problem. On the other hand, as seen in Theorem~\ref{theo:Sinkhorn}, the optimal transport plan of the regularized OT problem $\mathrm{OT}_{\gamma}$ is a rescaling of the Gibbs kernel $G_{ij} = \exp\np{-\frac{c_{ij}}{\gamma}}$ and the Sinkhorn iterates involve this kernel as well. When $\gamma$ is too close to $0$, Sinkhorn's algorithm shows numerical instabilities. Thus, we refrain from using a single regularizing parameter for every OT problem involved in Equation~\eqref{NestedDistance} and we simply put one that seems big enough to avoid numerical issues, namely we set $\gamma = \frac{\max_{ij}c_{ij}}{30}$, which changes as the cost matrix is updated. Hence the regularizing parameters do not explicitly appear in the notation $\mathrm{END}$ of the Entropic regularization of the Nested Distance. For every time $t\in \ce{1:T-1}$ and every \emph{node} $x_{1:t} \in \X_{1:t}$, define its set of \emph{children} $x_{1:t}^+ := \left\{\tilde{x}_{1:t+1} \in \X_{1:t+1} \mid \tilde{x}_{1:t} = x_{1:t} \right\}$.

\begin{definition}[Entropic regularization of Nested Distance between scenario trees]
 \label{def:END}
 Let $X$ and $Y$ be two scenario trees. Given $r\geq 1$, and the metric $d\np{x,y} = \Vert x-y \Vert_r$ over $\R^T$, for every $t\in \ce{1,T}$, compute recursively backward in time functions $c_t : \X_{1:T} \times \Y_{1:T} \to \overline{\R}$ by
 \begin{equation}
  \label{EntropicNestedDistance}
  \left\{
  \begin{aligned}
    & c_T\np{x_{1:T}, y_{1:T}} = d\np{x_{1:T}, y_{1:T}}, \ \forall \np{x_{1:T}, y_{1:T}} \in \X_{1:T} \times \Y_{1:T},                                                      \\
    & c_{t}\np{x_{1:T}, y_{1:T}} = \mathrm{OT}_{\gamma}\np{P_{t+1}\np{ \cdot \mid X_{1:t} = x_{1:t}}, \tilde{P}_{t+1}\np{\cdot \mid Y_{1:t} = y_{1:t}}; c_{t+1}^r } ^{1/r}, \\
    & \forall t\in \ce{1,T-1}, \ \forall \np{x_{1:T}, y_{1:T}} \in \X_{1:T} \times \Y_{1:T}, \ \gamma = \max_{\substack{x_{1:t+1} \in x_{1:t}^+                             \\ y_{1:t+1} \in y_{1:t}^+}} c_{t+1}^r\np{x_{1:{t+1}}, y_{1:t+1}} / 30.
  \end{aligned}
  \right.
 \end{equation}
 Set $\mathrm{END}_r\np{X, Y} := \mathrm{OT}_{\gamma}\np{P_T, \tilde{P}_T, c_1^r}^{1/r}$, with $\gamma = \max_{\substack{x_{1:1} \in \X_{1:1} \\ y_{1:1} \in \Y_{1:1}}} c_{t+1}^r\np{x_{1:{1}}, y_{1:1}} / 30$, it is the \emph{Entropic regularization of the $r$-Nested Distance} between the scenario trees $X$ and $Y$.
\end{definition}

As for the Nested Distance, the value of $c_t$ only depends of the nodes at time $t\leq T$. Moreover, note that by Remark~\ref{rem:ot_vs_regot}, for every $r\geq 1$ and scenario trees $X$ and $Y$,
\(
\mathrm{ND}_r\np{X, Y} \leq \mathrm{END}_r\np{X, Y}.
\)
Hence even though $\mathrm{END}_r$ is not a distance between scenario trees, it still quantifies proximity between scenario trees and maintain the main desirable feature of the Nested Distance: denoting by $v$ the value of a MSP satisfying the regularity assumptions of \cite[Theorem 11]{Pf.Pi2012}, there exists a constant $L>0$ such that for every scenario trees $X$ and $Y$ we have
\[
 \lvert v(X) - v(Y) \rvert \leq L\cdot \mathrm{ND}_r\np{X,Y} \leq L \cdot \mathrm{END}_r\np{X,Y}.
\]

\section{Numerical experiment}\
\label{sec:num}

We compare an implementation of the Nested Distance and an implementation of its regularized counterpart when $\R^N$, $N\in \N$, is endowed with the euclidean distance.

First, we randomly generate a scenario tree of given depth $T$ by a forward procedure. Starting from a root note, at each time step draw a uniformly random number of children between $1$ and a given number, here $3$. Every node at time $t$ has the given number of children whose values are random as well. The tree generation is done using the Julia package ScenTrees.jl, see \cite{Ki.Pi.Pf2020}. The discrete optimal transport problems are solved using the Julia package OptimalTransport.jl which uses a primal-dual interior point method. Second, we compute the Nested Distance and the Entropic regularization of the Nested Distance for pairs of tree generated as above. In Figure~\ref{result_num} we give the average of $10$ pairs of comparisons for a given horizon $T$.

We do a first batch (see Figure~\ref{result_num}) of simulations for the case of the entropic regularization of the $2$-Nested Distance, that is $r=2$. In Figure~\ref{result_num} the column "Relative error" represents the ratio $\frac{\mathrm{ND}_2 - \mathrm{END}_2}{\mathrm{END}_2}$. The results of Figure~\ref{result_num} show that even without carefully tuning the regularizing parameter $\gamma$ involved in each intermediate optimal transport problem, the Entropic regularization of the Nested Distance gives values that are close to its unregularized counterpart with an interesting speedup. We comfort this trend on a second batch (see Figure~\ref{result_num_2}) of simulations for the case $r=1$, however the relative error is degraded compared to the case $r=2$.

We refer the reader to the interactive (Julia Jupyter notebook) example available at \url{https://github.com/BenoitTran/END}. Note that some numerical instabilities arise for small values of $\gamma$, which could have been dealt with the log-sum-exp trick.

Tuning carefully the parameter $\gamma$ could yield better tradeoffs between speedup and precision. Moreover, a greedy variant of Sinkhorn's algorithm, called Greenhorn, could lead to even better and stable numerical results as shown in the works \cite{Al.Ni.Ri2017,pmlr-v97-lin19a}.

\begin{table}
 \centering
 \begin{tabular}{|c|c|c|c|c|}
  \hline
  Horizon $T$ & Time $\mathrm{ND}_2$ (ms) & Time $\mathrm{END}_2$ (ms) & Speedup & Relative error (\%) \\
  \hline
  2           & 0.26                      & 0.014                      & 16      & 0.14                \\
  \hline
  4           & 3.8                       & 0.14                       & 25      & 0.25                \\
  \hline
  6           & 115                       & 6.3                        & 33      & 0.51                \\
  \hline
  8           & 1077                      & 28                         & 35      & 0.35                \\
  \hline
  10          & 18205                     & 493                        & 36      & 0.41                \\
  \hline
 \end{tabular}
 \caption{\label{result_num} Average results after $10$ runs when $r=2$ with varying horizon $T$ and given maximum number of children of each node set to $3$. }
\end{table}

\begin{table}
 \centering
 \begin{tabular}{|c|c|c|c|c|}
  \hline
  Horizon $T$ & Time $\mathrm{ND}_1$ (ms) & Time $\mathrm{END}_1$ (ms) & Speedup & Relative error (\%) \\
  \hline
  2           & 0.14                      & 0.0075                     & 12.9    & 0.98                \\
  \hline
  4           & 4.6                       & 0.072                      & 31.7    & 2.92                \\
  \hline
  6           & 62                        & 0.92                       & 68.8    & 3.80                \\
  \hline
  8           & 2178                      & 27                         & 70.5    & 6.93                \\
  \hline
  10          & 4043                      & 94                         & 37.40   & 6.15                \\
  \hline
 \end{tabular}
 \caption{\label{result_num_2} Average results after $10$ runs when $r=1$ with varying horizon $T$ and given maximum number of children of each node set to $3$.}
\end{table}

\bibliographystyle{spmpsci}      
\bibliography{EntropicNestedDistance}   

\end{document}